\documentclass{amsart}          % Nicer than default article style:  less
\usepackage{amsmath,amsthm}     % Handy math stuff, theorem environments.
\usepackage{amssymb}            % Fancy math symbols.
\usepackage{euscript}           % Nice script font.
\usepackage{enumerate,calc}
\usepackage[matrix,arrow,curve,frame]{xy}    % XY-pic diagram pac

\xymatrixcolsep{1.9pc}                          % Adjust size of diagrams.
\xymatrixrowsep{1.9pc}
\newdir{ >}{{}*!/-9pt/\dir{>}}                  % Make better tailed arrows

\newtheorem{thm}[subsection]{Theorem}
\newtheorem{defn}[subsection]{Definition}
\newtheorem{prop}[subsection]{Proposition}

\newtheorem{cor}[subsection]{Corollary}
\newtheorem{lemma}[subsection]{Lemma}

\theoremstyle{definition}  % Bold headings and Roman body text.

\newtheorem{convention}[subsection]{Convention}

  {\end{list}}

\newcommand{\dfn}{\textbf} % Make defined words bold.

\newcommand{\mdfn}[1]{\dfn{\mathversion{bold}#1}} % Even make math bold

\newcommand{\Wedge}             {\vee}

\newcommand{\tens}              {\otimes}               %tensor
  
\newcommand{\cat}{\EuScript}    % Use \EuScript to name a category.
\newcommand{\cC}{{\cat C}}

\newcommand{\cR}{{\cat R}}
\newcommand{\cS}{{\cat S}}

\newcommand{\sSet}{s{\cat Set}}

\DeclareMathOperator*{\colim}{colim}

\DeclareMathOperator{\Hom}{Hom}

\DeclareMathOperator{\Ex}{Ex}
\DeclareMathOperator{\sd}{sd}
\DeclareMathOperator{\Cone}{Cone}
\DeclareMathOperator{\im}{im}

\newcommand{\ra}{\rightarrow}                   % right arrow
\newcommand{\lra}{\longrightarrow}              % long right arrow
\newcommand{\lla}{\longleftarrow}               % long left arrow
\newcommand{\llla}[1]{\stackrel{#1}{\lla}}      % labeled long right
\newcommand{\inc}{\hookrightarrow}              % inclusion
\newcommand{\ovcat}{\downarrow}

\newcommand{\restr}[1]{\!\mid_{#1}}

\newcommand{\bd}[1]{\partial\Delta^{#1}}

\newcommand{\Lamb}[1]{\Lambda^{#1}}
\newcommand{\del}[1]{\Delta^{#1}}

\newcommand{\rea}[1]{|{#1}|}             %geometric realization of #1
\newcommand{\map}{\rightarrow}

\newcommand{\ceck}[1]{\Cech(#1)}         %Cech complex for #1
\newcommand{\oceck}[1]{\Cech^{o}(#1)}    %Ordered Cech complex for #1
\newcommand{\oreal}[1]{\rea{\oceck{U}}}  %Realization of ordered Cech cplex
\newcommand{\creal}[1]{\rea{\ceck{U}}}   %Realization of the Cech complex
\newcommand{\Cech}{\check{C}}

\newcommand{\id}{\text{id}}

\numberwithin{equation}{subsection}

\newenvironment{myequation}
  {\addtocounter{subsection}{1}\begin{eqnarray}}
  {\end{eqnarray}$\!\!$}

\begin{document}

\title{Weak equivalences of simplicial presheaves}

\author{Daniel Dugger} 
\author{Daniel C. Isaksen}

\address{Department of Mathematics\\ Purdue University\\ West
Lafayette, IN 47907 } 

\address{Department of Mathematics\\ University of Notre Dame\\
Notre Dame, IN 46556}

\email{ddugger@math.purdue.edu}

\email{isaksen.1@nd.edu}

\begin{abstract}
Weak equivalences of simplicial presheaves are usually defined in
terms of sheaves of homotopy groups.  We give another characterization
using relative-homotopy-liftings, and develop the tools necessary to
prove that this agrees with the usual definition.  From our lifting
criteria we are able to prove some foundational (but new) results about
the local homotopy theory of simplicial presheaves.
\end{abstract}

\maketitle

%\tableofcontents

\section{Introduction}
In developing the homotopy theory of simplicial sheaves or presheaves,
the usual way to define weak equivalences is to require that a map
induce isomorphisms on all sheaves of homotopy groups.  This is a
natural generalization of the situation for topological spaces, but
the `sheaves of homotopy groups' machinery (see
Definition~\ref{de:Ill-we}) can feel like a bit of a mouthful.  The
purpose of this paper is to unravel this definition, giving a fairly
concrete characterization in terms of lifting properties---the kind of
thing which feels more familiar and comfortable to the ingenuous
homotopy theorist.

The original idea came to us via a passing remark of Jeff Smith's: He
pointed out that a map of spaces $X \ra Y$ induces an isomorphism on
homotopy groups if and only if every diagram
\begin{myequation}
\label{eq:we}
 \xymatrixcolsep{1.5pc}\xymatrix{
  &S^{n-1} \ar[r]\ar@{ >->}[d]\ar@{ >->}[dl] & X \ar[d] \\
  D^n \ar@{ >->}[d] \ar@{.>}[urr]  
             & D^n \ar[r]\ar@{ >->}[dl] & Y  \\
  D^{n+1}   \ar@{.>}[urr]
}
\end{myequation}
admits liftings as shown (for every $n\geq 0$, where by convention we
set $S^{-1}=\emptyset$).  Here the maps $S^{n-1} \inc D^n$ are both
the boundary inclusion, whereas the two maps $D^n \inc D^{n+1}$ in the
diagram are the two inclusions of the surface hemispheres of
$D^{n+1}$.  The map $D^{n+1}\ra Y$ should be thought of as giving a
homotopy between the two maps $D^n \ra Y$ relative to $S^{n-1}$.  In
essence, the above lifting condition just guarantees the vanishing of
the relative homotopy groups of $X\ra Y$.

One advantage of this formulation is that one doesn't have to worry
about basepoints, but it also has other conveniences.  If one looks
back on the classical lifting theorems in \cite{Sp}, for instance, it
is really the above property---rather than the isomorphism on homotopy
groups---which is being made use of over and over again.  In working
with simplicial presheaves, it eventually became clear that a version
of the above characterization was a useful thing to have around.  It
comes in at several points in \cite{DHI}, where it is used to
inductively produce liftings much like in \cite{Sp}.

Whereas for topological spaces the above characterization is
`obvious', for simplicial presheaves it requires a little bit of work.
Intuitively the result is clear, but to actually write down a proof
one must (a) struggle with the combinatorics of simplicial sets, and
(b) deal with the `local homotopy theory' which demands that
everything be accomplished by a finite number of lifting arguments.
The trouble is that the modern way of avoiding (a) is to use the tools
of model categories, but because of (b) we don't have these at our
disposal.  In this paper we develop some basic machinery for handling
this situation, so that in the end one can write out the proof fairly
smoothly.

\medskip

To describe the results more explicitly, we'll first give an analagous
characterization for weak equivalences $X\ra Y$ of simplicial sets.
In this case, we must assume that $X$ and $Y$ are fibrant.  The
$(n-1)$-sphere is replaced by $\bd{n}$; the $n$-disks are replaced by
$\del{n}$; and $D^{n+1}$ is replaced by the pushout $RH(\del{n},
\bd{n})$ of the diagram
\[
\bd{n} \leftarrow \bd{n} \times \del{1} \map \del{n} \times \del{1}.
\]
The simplicial set $RH(\del{n}, \bd{n})$ is the domain of simplicial
homotopies---relative to $\bd{n}$---between maps out of $\del{n}$.
Once these substitutions are made into diagram (\ref{eq:we}), one gets
the same criterion for the map $X\ra Y$ to be a weak equivalence of
simplicial sets.  See Proposition~\ref{pr:ss-we}.

The generalization to the case of simplicial presheaves is now
reasonably transparent.  The main result of the paper is the
following:

\begin{thm}
\label{th:intromain}
Let $\cC$ be a Grothendieck site, and let $F\ra G$ be a map between
locally fibrant simplicial presheaves.  Then $F\ra G$ induces an
isomorphism on all sheaves of homotopy groups (for all choices of
basepoint) if and only if it has the following property: for every
solid-arrow diagram
\begin{myequation}
\label{eq:pre-we}
 \xymatrixcolsep{1.5pc}\xymatrix{
  &\bd{n} \tens X \ar[r]\ar@{ >->}[d]\ar@{ >->}[dl] & F \ar[d] \\
  \del{n} \tens X \ar@{ >->}[d] \ar@{.>}[urr]  
             & \del{n} \tens X \ar[r]\ar@{ >->}[dl] & G  \\
  RH(\del{n}, \bd{n}) \tens X   \ar@{.>}[urr]
}
\end{myequation}
in which $X$ is representable, there exists a covering sieve of $X$
such that for any $U \map X$ in the sieve, 
the diagram obtained by pulling back to $U$ has liftings as shown.
\end{thm}
(The dotted arrows in the above diagram are called
`relative-homotopy-liftings', and the fact that they only exist
locally leads us to call this property the `local RHLP'---see
Defintion~\ref{de:llp}, as well as the general discussion in
Section~\ref{se:genprops}).

The advantage of this viewpoint on weak equivalences is that it
provides a framework for using lifting arguments instead of
computations of homotopy groups.  For questions in homotopy theory,
lifting is sometimes a more convenient tool.  In our case, this
lifting characterization will be used to give elegant proofs of
various results (some old, some new) about local weak equivalences and
local-fibrations.  These are given in Section~\ref{se:app}. To readers
of \cite{J1,J3}, these come as no surprise.  However, our proofs seem
simpler and more conceptual than the ones involving sheaves of
homotopy groups, in particular avoiding all references to stalks or
Boolean localizations.  And in some cases we don't know any proof other
than via the lifting criterion.

We should remark that this approach via liftings is not at all meant
to replace the definition involving sheaves of homotopy groups---in
some situations that is exactly the tool that is needed.  But in
general it is good to have both descriptions at one's disposal.

The chief motivation for writing this paper was its application to our
study of localization for simplicial presheaves \cite{DHI}.  But we've
also found that the techniques of homotopy-liftings are convenient
tools to have around, and should be more well-known among abstract
homotopy theorists.  Reedy \cite{R} worked with a dual version of
these lifting criterion in the context of abstract model categories,
and used them to prove several key lemmas.  We have reproduced a
couple of his proofs here for completeness, and with the goal of
popularizing these ideas.

\subsection{Organization of the paper} \mbox{}

Sections \ref{se:backgrnd}--\ref{se:combprf} deal only with simplicial
sets.  Section 2 has a few background results, and then in Section
\ref{se:genprops} we define relative-homotopy-liftings and develop
their basic properties.  In Section \ref{se:we-ss} these ideas are
applied to get a lifting criterion for weak equivalences of simplicial
sets.  Unfortunately, one of the key steps is easiest to prove using
model-category theoretic methods, and these do not generalize to the
simplicial presheaf setting.  So Section~\ref{se:combprf} is devoted
to giving a completely combinatorial proof for this result.  We
finally get to simplicial presheaves in Section~\ref{se:we}.  We
recall the traditional definition of weak equivalence using sheaves of
homotopy groups, and then prove Theorem~\ref{th:intromain}.  Section
\ref{se:app} concludes with some applications of this theorem.

We assume that the reader is familiar with standard results from the
homotopy theory of simplicial sets, including material to be found in
\cite{M} or \cite{GJ}.  Much of what we discuss, especially in
Sections \ref{se:genprops} and \ref{se:we-ss}, can be easily
generalized to abstract model categories, but we do not treat this
extra generality here.  We'll also assume a familiarity with sheaf
theory and the homotopy theory of simplicial presheaves, for which we
refer the reader to \cite{J1}.  It should be clear from our arguments
how indebted we are to that paper.

\section{Background on simplicial sets} 
\label{se:backgrnd}

We start with some basic facts.  Let $S^n$ be the sphere $\del{n} /
\bd{n}$.  If $(K,x)$ is a pointed simplicial set, then $\pi_n(K,x)$
denotes the set of maps $(\del{n},\bd{n})\ra (K,x)$ modulo the
equivalence relation generated by simplicial homotopy relative to
$\bd{n}$.  Of course this set has homotopical meaning only if $K$ is
fibrant.

In \cite{K}, Kan constructed a fibrant-replacement functor called
$\Ex^\infty$.  First, let $\sd$ be the barycentric subdivision functor
\cite[p.~183]{GJ}.  For any simplicial set $X$, $\Ex X$ is the simplicial set
whose $k$-simplices are elements of the set $\Hom(\sd \del{k}, X)$.
The functor $\Ex$ is right adjoint to $\sd$.  Now $\Ex^n$ is the
$n$-fold composition of $\Ex$, and $\Ex^\infty$ is $\colim_n \Ex^n$.
The functor $\Ex^\infty$ has some nice properties one wouldn't expect
from an arbitrary fibrant-replacement functor: It preserves
fibre-products, it preserves the set of $0$-simplices, and it
preserves fibrations.  These properties all follow immediately from
the definition.

The following two basic lemmas about simplicial sets will be used later.
The first is obvious, but
its statement and proof becomes important when considering simplicial
presheaves later.  The point is that the proof uses 
only basic lifting properties, not fancy model theoretic results.

\begin{lemma}
\label{le:contract-extn}
Let $i\colon K \inc L$ be a cofibration, and let $X$ be a fibrant
simplicial set.
If $f\colon K \map X$ factors through any contractible simplicial set $M$,
then $f$ is simplicially null-homotopic and $f$ extends over $i$.
\end{lemma}
% True when $X$ is a locally-fibrant presheaf as long as $K$ and $L$ are
% finite.

\begin{proof}
For any simplicial set $Y$, let $\Cone Y$ be 
$(Y \times \del{1}) / (Y \times \{1\})$.  
We have a diagram
\[
\xymatrix{
K \ar[r] \ar[d] & M \ar[d] \ar[r] & X \\
\Cone K \ar[r] & \Cone M. \ar@{.>}[ur]                          }
\]
The map $M \map \Cone M$ is an acyclic cofibration because $M$ is
contractible, so there is a lift as shown.  Composition with 
$\Cone K \map \Cone M$ gives the desired simplicial null-homotopy.

Now we have a diagram
\[
\xymatrix{
K \ar[r] \ar[d] & \Cone K \ar[d] \ar[r] & X \\
L \ar[r] & \Cone L. \ar@{.>}[ur]                          }
\]
The map $\Cone K \map \Cone L$ is an acyclic cofibration,
so there is a lift as shown.  Composition with 
$L \map \Cone L$ gives the desired extension over $i$.
\end{proof}

Note that if $K$, $L$, and $M$ are finite simplicial sets, then the
desired lift can be produced using only finitely many applications of
the Kan extension condition for $X$.  This will be important when we
start generalizing to simplicial presheaves.

\begin{lemma}
\label{le:a-cofib}
Let $i: K \map L$ be an acyclic cofibration between finite simplicial
sets.  Then $i$ can be built from the maps $\Lambda^{n,k} \inc \del{n}$ 
by a finite
number of retracts, cobase changes, and compositions.
\end{lemma}

\begin{proof}
We know that $i$ is a retract of a relative $J$-cell complex $j: M
\map N$ (cf. \cite[Def. 12.5.8]{H}), where $J$ is the set of maps of
the form $\Lambda^{n,k} \inc \del{n}$.  Since $L$ is finite, its image
in $N$ belongs to a finite subcomplex.  Thus $i$ is actually a retract
of a finite relative $J$-cell complex.
\end{proof}

%%%%%%%%%%%%%%%%%%%%%%%%%%%%%%%%%%%%%%%%%%%%%%%%%%%%%%%%%%%%%%%%%

\section{Generalities about homotopy-liftings}
\label{se:genprops}

This section establishes the definition and basic properties of what
we call `relative-homotopy-liftings'.  

\begin{defn}
A square of simplicial sets
\begin{myequation}
\label{eq:rhlp}
\xymatrix{
K \ar[r] \ar[d] & X \ar[d] \\
L \ar[r] & Y            }
\end{myequation}
is said to have a \mdfn{relative-homotopy-lifting}
if there exists a map $L \map X$ such that the upper left triangle commutes
and there is a simplicial homotopy relative to $K$ from the composition
$L \map X \map Y$ to the given map $L \map Y$.

The map $X\ra Y$ has the \dfn{relative-homotopy-lifting property}
(RHLP) with respect to $K\ra L$ if every square (\ref{eq:rhlp}) has a
relative-homotopy-lifting.
\end{defn}

Reedy \cite[Lem.~2.1]{R} used the dual to the above definition.  Like
him, we could haved defined this property in an arbitrary simplicial
model category (one probably doesn't even need the model category to
be simplicial).  All of our basic results go through in that
generality, but we won't ever need this.

Given $K\ra L$, let \mdfn{$RH(L,K)$} denote the pushout of the diagram
\[        K   \llla{\pi}   K\times\del{1} \lra L\times\del{1}
\]
where the left map is the projection.
The notation stands for `Relative-Homotopy':
To give a map $RH(L,K)\ra X$
means precisely to give two maps $L \ra X$ which agree on $K$,
together with a simplicial homotopy between them relative to $K$.  Note that
there is a canonical map $L\amalg_K L \ra RH(L,K)$.  This map is
a cobase change of the map 
\[
(K \times \del{1}) \cup (L \times \bd{1}) \map L \times \del{1},
\]
so it is a cofibration if $K\ra L$ is.  

We will sometimes use the fact that the existence of
relative-homotopy-liftings can be rephrased as
saying that the diagram
\begin{myequation}
\label{eq:catrhlp}
 \xymatrixcolsep{1.5pc}\xymatrix{
  &K \ar[r]\ar[d]\ar[dl] & X \ar[d] \\
  L \ar[d] \ar@{.>}[urr]  
             & L \ar[r]\ar[dl] & Y  \\
  RH(L,K)   \ar@{.>}[urr] }
\end{myequation}
admits liftings as shown.  While this diagram may seem somewhat
awkward (especially when seeing it for the first time), it is often a
very useful tool.

Here are some basic properties of relative-homotopy-liftings:

\begin{lemma}
\label{le:hlift-comp}
Let $f\colon X \map Y$ be a fixed map of simplicial sets.  Consider
the class of all maps $K \map L$ with respect to which $f$ has the
RHLP.  This class is closed under cobase changes and retracts.  If $Y$
is fibrant, then the composition of two cofibrations in the class is
still in the class.
\end{lemma}

\begin{proof}
Closure under cobase changes follows from consideration of
(\ref{eq:catrhlp}) and the fact that $RH(L \amalg_K M, M)$ is
isomorphic to $RH(L, K) \amalg_K M$.  Closure under retracts follows
from the usual formal argument with lifting properties.

For composition, we start
with two cofibrations $i\colon K \map L$ and $j\colon L \map M$ in the
class.  Consider a lifting problem
\[ \xymatrix{ K \ar[r]^g\ar[d]_{ji} & X \ar[d]^{f} \\
              M \ar[r]_{h} & Y.
}
\]
The first step is to produce a homotopy-lifting
\[ \xymatrix{ K \ar[r]^{g}\ar[d]_{i} & X \ar[d]^{f} \\
              L \ar[ur]^{l}\ar[r]_{hj} & Y
}
\]
relative to $K$.
Let $H:RH(L, K) \map Y$ be the relative-homotopy from
$fl$ to $hj$.
So now we look at the diagram 
\[ \xymatrix{
   (M\times\{1\}) \amalg_{L \times \{1\}} RH(L,K)
             \ar[r]^-{h \amalg H}  \ar@{ >->}[d]_{\sim} 
    & Y\\
    RH(M,K) \ar@{.>}[ur]_J
}
\]
and produce a lifting $J$ using that $Y$ is fibrant.
Note that the vertical map above is an acyclic cofibration because
it is a cobase change of the acyclic cofibration
\[
(M\times\{1\}) \cup (L\times \del{1}) \map M \times \del{1}.
\]
Let $m$ be the map $J|_{M \times \{0\}}$.  
Note that $m$ is simplicially homotopic
to $h$ relative to $K$.

At this point we have the square
\[ \xymatrix{ L \ar[r]^l\ar[d]_{j} & X\ar[d]^{f} \\
              M \ar[r]^m & Y.
}
\]
We produce a relative-homotopy-lifting $n$.
A diagram chase shows that $nji$ equals $g$.  On the other hand,
$fn$ is simplicially homotopic
to $m$ relative to $K$ and hence also to $h$; here we use that
relative-homotopy is transitive because $Y$ is fibrant.
\end{proof}

\begin{cor}
\label{co:hlift-comp}
If $f\colon X\ra Y$ has the RHLP with respect to the maps $\bd{n}\inc
\del{n}$ for all $n\geq 0$ and $Y$ is fibrant, 
then $p$ also has the RHLP with respect to
all cofibrations $K\inc L$ of finite simplicial sets.
\end{cor}

\begin{proof}
Every such cofibration $K \inc L$ can be constructed by a finite number
of compositions and cobase changes from the generating
cofibrations $\bd{n} \inc \del{n}$.
\end{proof}

\begin{prop}
\label{pr:homotop}
Suppose that $K\inc L$ is a cofibration and we are
given a square $\cS$ of the form
\[
\xymatrix{
K\times\del{1} \ar[r]^-{H_K} \ar[d] & X \ar[d]^f \\
L\times\del{1} \ar[r]^-{H_L} & Y            }
\]
in which $X$ and $Y$ are fibrant.
Let $\cS_0$ denote the square obtained by restricting $H_K$ and $H_L$
to time $t=0$, and similarly for $\cS_1$.
Then $\cS_0$ has a relative-homotopy-lifting if and only if $\cS_1$ does.
\end{prop}

\begin{proof}
Suppose that $\cS_0$ has a relative-homotopy-lifting $l_0$.  Then by gluing
$l_0$ to $H_K$ we get
$(L\times \{0\}) \cup (K\times\del{1}) \ra X$, and since $X$ is
fibrant this map extends over $L\times\del{1}$.  Let $l_1$ denote the
restriction of this map to $L\times \{1\}$; we will show that $l_1$ is
the desired relative-homotopy-lifting.

Pushing the homotopy $L\times\del{1}\ra X$ down into $Y$, we can
glue it to $H_L$ to get a map
defined on $(L\amalg_K L)\times\del{1}$.  Together with the 
relative-homotopy from $fl_0$ to $H_L|_{t=0}$, 
we find that we actually have a map
$[RH(L,K)\times\{0\}] \cup [(L\amalg_K L)\times\del{1}] \ra Y$.  Since
$Y$ is fibrant, this extends over $RH(L,K) \times \del{1}$.
Restricting to time $t=1$ gives the desired relative-homotopy of
$fl_1$ with $H_L|_{t=1}$.
\end{proof}

\begin{prop}
\label{pr:fibhl}
Let $f\colon X\ra Y$ be a map between fibrant simplicial sets.  Then
$f$ has the RHLP with respect to 
every acyclic cofibration $K \inc L$.
In particular, $f$ has the RHLP with respect to
the maps $\Lambda^{n,k} \inc \del{n}$.
\end{prop}
% Generalizes to locally-fibrant presheaves when $A$ and $B$ are finite.

\begin{proof}
Consider a square
\[
\xymatrix{
K \ar[r]^g \ar[d] & X \ar[d]^f \\
L \ar[r]_h & Y.                 }
\]
Using that $X$ is fibrant,
there is a map $l\colon L\ra X$ extending $K\ra X$.  
We must give a relative-homotopy from $fl$ to $h$.

The cofibration $L \amalg_K L \inc RH(L, K)$ is an acyclic
cofibration because it is a cobase change of the acyclic cofibration
$(L\times \bd{1}) \cup (K\times\del{1}) \inc L\times\del{1}$.  
We may extend the map $fl \amalg h: L \amalg_K L \map Y$ to
$RH(L, K) \ra Y$ because $Y$ is fibrant.  This gives us the desired
relative-homotopy.
\end{proof}

Observe that if $K \ra L$ is $\Lambda^{n,k}\ra\del{n}$ then the
proof only requires a finite number of uses of the Kan extension
condition.  Combined with Lemma~\ref{le:a-cofib}, this tells us
that the proposition applies to the simplicial presheaf setting
whenever $K$ and $L$ are finite simplicial sets.  This won't be needed
until Section~\ref{se:we}.

%%%%%%%%%%%%%%%%%%%%%%%%%%%%%%%%%%%%%%%%%%%%%%%%%%%%%%%%%%%%%%%%%%%%%

\section{Weak equivalences of simplicial sets}
\label{se:we-ss}

The following proposition now shows that weak equivalences between
fibrant simplicial sets can be detected using relative-homotopy-liftings.  
This is analogous to the situation discussed in the introduction
for topological spaces, 
where every object is fibrant.

\begin{prop}
\label{pr:ss-we}
A map $f\colon X \map Y$ between fibrant simplicial sets 
is a weak equivalence if and only if 
it has the RHLP with respect to the maps $\bd{n}\inc\del{n}$, for all
$n\geq 0$.
\end{prop}

When the simplicial sets are not fibrant one has to allow oneself to
subdivide $\bd{n}$ and $\Delta^n$, but we won't pursue this.

The following proof is similar to the proof of \cite[Lem.~2.1]{R}.  
The difference is that we only consider the RHLP with
respect to the generating cofibrations, while Reedy considers
the RHLP with respect to all cofibrations.
We include the full details for completeness.

\begin{proof}
First suppose that $f$ has the RHLP.  By Corollary
\ref{co:hlift-comp}, $f$ has the RHLP with respect to the cofibrations
$* \inc S^n$ for all $n \geq 1$,
as well as $\emptyset\ra *$.  
This shows that $\pi_n X \map \pi_n Y$ is surjective (for any choice of
basepoint).  Similarly, $f$ has the RHLP with respect to the
cofibrations $S^n \Wedge S^n \inc RH(S^n, *)$ for all $n$.  This shows
that $\pi_n X \map \pi_n Y$ is injective.

Conversely, we'll now suppose that $f$ is a weak equivalence.
Consider a square
\[
\xymatrix{
\bd{n} \ar[r] \ar@{ >->}[d] & X \ar[d] \\
\Delta^{n} \ar[r] & Y.            }
\]
Factor $f$ into an acyclic cofibration $i\colon X \inc Z$
followed by an acyclic fibration $p \colon Z \map Y$.  Since $X$ is
fibrant, there is a map $g\colon Z \map X$ making $X$ a retract of
$Z$.  Now the square
\[
\xymatrix{
\bd{n} \ar[r] \ar@{ >->}[d] & X \ar[r]^i &  Z\ar[d]^p \\
\Delta^{n} \ar[rr] & & Y            }
\]
has a lift $h$ because $p$ is an acyclic fibration.

The composition $gh$ is the desired homotopy-lift.  Using that $gi$ is
the identity, the upper left triangle commutes.  Working in the
undercategory $\bd{n}\ovcat \sSet$, we see that $ig$ represents the
same map as $\id_Z$ in the homotopy category---therefore $pigh$
represents the same map as $ph$ in the homotopy category.  But these
latter two maps have cofibrant domain and fibrant target, so they
are actually simplicially homotopic in $\bd{n}\ovcat \sSet$.  The
simplicial set $RH(\del{n},\bd{n})$ is precisely a cylinder object for
$\del{n}$ in this undercategory, so $pigh$ and $ph$ are
simplicially homotopic relative to $\bd{n}$.
\end{proof}

Reedy \cite[Th.~B]{R} showed that base changes along fibrations
preserve weak equivalences between fibrant objects.  His proof used
the criterion of Proposition \ref{pr:ss-we} (suitably generalized to
arbitrary model categories) to detect weak equivalences.  We use this
idea to obtain the following elementary proof of right properness for
simplicial sets---most standard references \cite{GJ,H} use topological
spaces to prove this.
We include full details because this same proof will be applied to the
case of simplicial presheaves.

\begin{cor}[Right properness]
\label{co:ss-rp}
Let $f:X \map Y$ be a weak equivalence of simplicial sets, and let
$p:Z \map Y$ be a fibration.  Then the map $X \times_Y Z \map Z$
is also a weak equivalence.
\end{cor}

\begin{proof}
We need only show that $\Ex^\infty (X \times_Y Z) \map \Ex^\infty Z$
is a weak equivalence.
Note that $\Ex^\infty$ commutes with fibre-products
and preserves fibrations, so this map is a base change of
the weak equivalence $\Ex^\infty f$ along the fibration
$\Ex^\infty p$.
Therefore,
we may assume that $X$, $Y$ and $Z$ are already fibrant.

The rest of the proof is the same as Reedy's argument.
Suppose given a square
\[
\xymatrix{
\bd{n} \ar[r]^-g \ar@{ >->}[d] & X \times_Y Z \ar[d] \\
\del{n} \ar[r]_h & Z.                        }
\]
We want to find a relative-homotopy-lifting for this square.
First, take a relative-homotopy-lifting $l$ for the composite square
\[
\xymatrix{
\bd{n} \ar[r]^-g \ar@{ >->}[d] & X \times_Y Z \ar[r] & X \ar[d]^f \\
\del{n} \ar[r]_h & Z \ar[r]_p & Y,                        }
\]
which exists by Proposition \ref{pr:ss-we}
because $X \map Y$ is a weak equivalence between fibrant
simplicial sets.
Now consider the square
\[
\xymatrix{
\del{n} \times \{1\} \ar[r]^-h \ar@{ >->}[d] & Z \ar[d] \\
RH(\del{n}, \bd{n}) \ar[r] & Y,               }
\]
where the bottom horizontal map is the relative-homotopy from $fl$ to $ph$.
This square has a lift $H$ because the left vertical arrow is an acyclic 
cofibration.
Let $H_0$ and $H_1$ be the restrictions of $H$ to $\del{n} \times \{0\}$
and $\del{n} \times \{1\}$ respectively.  Note that $H_1 = h$ and
$pH_0 = fl$.

The maps $l$ and $H_0$ together define a map $m:\del{n} \map X
\times_Y Z$.  A diagram chase shows that $g$ is the restriction of $m$
to $\bd{n}$, and $H$ is the necessary relative-homotopy.
\end{proof}

%%%%%%%%%%%%%%%%%%%%%%%%%%%%%%%%%%%%%%%%%%%%%%%%%%%%%%%%%%%%%%%%%%%%

\section{A combinatorial proof}
\label{se:combprf}

In the previous section, Proposition \ref{pr:ss-we} compared weak
equivalences between fibrant simplicial sets to maps that have the
RHLP with respect to the cofibrations $\bd{n} \inc \del{n}$.
Unfortunately, the proof of one implication of the proposition relied
on model category theoretic methods.  When we generalize to simplicial
presheaves later on, these methods are not at our disposal.  Thus, our
goal in this section is to show by purely combinatorial methods that
if $f$ is a weak equivalence then it has the RHLP with respect to the
maps $\bd{n}\inc\del{n}$.

Throughout this section $f\colon X\ra Y$ denotes a map between
fibrant simplicial sets.

\medskip

First note that surjectivity on homotopy groups says precisely that
$f$ has the RHLP with respect to the maps $* \map S^{n}$ for all $n\geq 1$
as well as the map $\emptyset \ra *$.
Using this, we have:

\begin{lemma}
\label{le:lift1}
If $f\colon X \map Y$ is a map between fibrant simplicial sets that
induces surjections on homotopy groups, then $f$ has the RHLP with
respect to the maps $\Lambda^{n,k}\inc \bd{n}$, for any $n\geq 1$.
\end{lemma}

\begin{proof}
Suppose given a square
\[
\xymatrix{
\Lambda^{n,k} \ar[r] \ar@{ >->}[d] & X \ar[d] \\
\bd{n} \ar[r] & Y.                    }
\]
Since $\Lambda^{n,k}$ is contractible and $X$ is fibrant, the map
$\Lambda^{n,k} \ra X$ is simplicially null-homotopic (by
Lemma~\ref{le:contract-extn})---choose a null-homotopy.  By composing
with $X\ra Y$, we also get a null-homotopy for the composition
$\Lambda^{n,k} \map Y$; so we have a map $(\bd{n} \times \{0\}) \cup
(\Lambda^{n,k} \times \del{1}) \map Y$ that is constant on
$\Lambda^{n,k} \times \{1\}$.  Because $Y$ is fibrant, this map
extends to a map $\bd{n} \times \del{1} \map Y$ that is constant on
$\Lambda^{n,k} \times \{1\}$.

We have constructed a homotopy (in the sense of Proposition \ref{pr:homotop})
between the original square and a square of the form
\[
\xymatrix{
\Lambda^{n,k} \ar[r] \ar@{ >->}[d] & \mbox{}* \ar[r] \ar@{ >->}[d] 
    & X \ar[d] \\
\bd{n} \ar[r] & \bd{n} / \Lambda^{n,k} \ar[r] & Y.              }
\]
By Proposition \ref{pr:homotop}, we need only construct a 
relative-homotopy-lifting for this new square.  
The left square is a pushout, so 
we need only
construct a relative-homotopy-lifting for the right-hand square
by Lemma \ref{le:hlift-comp}.

Note that $\bd{n} / \Lambda^{n,k}$ is isomorphic to $S^{n-1}$.  Therefore,
$f$ has the RHLP with respect to $* \map \bd{n} / \Lambda^{n,k}$ 
because $f$ induces a surjection on $(n-1)$st homotopy groups.
\end{proof}

\begin{thm}
\label{th:comb-ss}
If $f\colon X\ra Y$ is a weak equivalence between fibrant
simplicial sets, then it has the RHLP with
respect to the maps $\bd{n}\inc\del{n}$, for all $n\geq 0$.
\end{thm}

\begin{proof}
Surjectivity on $\pi_0$ immediately gives the result for $n=0$.
So suppose $n\geq 1$ and we have a lifting diagram
\[\xymatrix{
\bd{n} \ar[r]^g \ar@{ >->}[d] & X \ar[d] \\
\del{n} \ar[r]^h & Y.}
\]
Routine lifting arguments show that there is a simplicial homotopy 
$\bd{n}\times \del{1} \ra X$ between $g$ and a map that
factors through $\bd{n}/\Lambda^{n,n}$.  
As in the proof of Lemma \ref{le:lift1}, we can extend this
to a simplicial homotopy $\del{n} \times \del{1} \map Y$, and
we are reduced to producing a
relative-homotopy-lifting for a square of the form
\[\xymatrix{
\bd{n}/\Lambda^{n,n} \ar[r]^-g \ar@{ >->}[d] & X \ar[d] \\
\del{n}/\Lambda^{n,n} \ar[r]^-h & Y.}
\]

Note that $\bd{n}/\Lambda^{n,n}$ is isomorphic to $S^{n-1}$ and
$\del{n}/\Lambda^{n,n}$ is contractible.  
Lemma \ref{le:contract-extn} 
shows that $g\colon S^{n-1}\ra X$ becomes null in $\pi_{n-1}(Y)$.
Since $f$ is injective on homotopy groups, $g$ is simplicially
null-homotopic.  Therefore, $g$ extends to a map
$l:\del{n}/\Lambda^{n,n} \map X$ by Lemma \ref{le:contract-extn}.
%In this case the $n$th face of $g$ represents a map
%$(\del{n-1},\bd{n-1}) \ra (X,*)$ which extends over
%$\del{n}/\Lambda^{n,n}$ when we compose with $p$.  By Lemma~\ref{le:inj},
%it extends over $\del{n}/\Lambda^{n,n}$ in $X$ as well: so we have a
%map $h\colon(\del{n},\Lambda^{n,n}) \ra (X,*)$ extending $g$.

Define a map $H\colon\bd{n+1} \ra Y$ by making the
$(n+1)$st face equal to $h$, the $n$th face equal to $fl$, and all the
other faces equal to the basepoint $*$.  Similarly, define $J
\colon \Lambda^{n+1,n+1} \ra X$ by making the $n$th face equal to $l$
and all the other faces equal to the basepoint.  So we have a 
square
\[
\xymatrix{ \Lambda^{n+1,n+1} \ar[r]^J \ar@{ >->}[d] & X \ar[d] \\
            \bd{n+1} \ar[r]_H & Y,
}
\]
and by Lemma~\ref{le:lift1} this has a
relative-homotopy-lifting $m$.  The $(n+1)$st face of $m$ gives a map
$\del{n} \ra X$ that is a relative-homotopy-lifting for our original
square.  
\end{proof}

%%%%%%%%%%%%%%%%%%%%%%%%%%%%%%%%%%%%%%%%%%%%%%%%%%%%%%%%%%%%%%%%%%%%
\section{Local weak equivalences of simplicial presheaves}
\label{se:we}

In this section we prove the main theorem (stated here as
Theorem~\ref{th:main}).  We start by recalling some of the tools from
\cite{J1}: the use of local lifting properties and sheaves of homotopy
groups.  Then we set up the local version of
relative-homotopy-liftings, and observe that everything we've done so
far still works in this setting.

\subsection{Local-liftings}

Fix a Grothendieck site $\cC$.  Recall that a map of simplicial
presheaves $F\ra G$ is a \dfn{local-fibration} if it has the following
property: given any square
\begin{myequation}
\label{di:locfib1}
\xymatrix{ \Lamb{n,k} \tens X \ar[r]\ar@{ >->}[d] & F \ar[d] \\
              \del{n} \tens X \ar[r] & G
}
\end{myequation}
in which $X$ is representable, 
there exists a covering sieve of $X$ such that for any  map $U\ra
X$ in the sieve, the induced diagram
\begin{myequation}
\label{di:locfib}
\xymatrix{ \Lamb{n,k}\tens U \ar[r] \ar@{ >->}[d] 
                 &\Lamb{n,k} \tens X \ar[r] & F \ar[d] \\
              \del{n}\tens U \ar[r]\ar@{.>}[urr]\ar[r]
                &\del{n} \tens X \ar[r] & G 
}
\end{myequation}
has a lifting as shown.  We are {\it not\/} requiring that
the liftings for different $U$'s be compatible in any way, only that
they exist.  
This kind of `local lifting property' will appear often
in the course of the paper, so we adopt the following convention:

\begin{convention}
Suppose given a lifting diagram like (\ref{di:locfib1}), in which a
representable presheaf $X$ appears.  We say this diagram has
\dfn{local liftings} if there exists a covering sieve $R$ of $X$ such
that for any $U\ra X$ in $R$ the diagram obtained by pulling back to
$U$ admits liftings.  For instance, using this language, a map $F\ra
G$ is a local-fibration provided that every diagram (\ref{di:locfib1})
admits local liftings.
\end{convention}

Because of Lemma \ref{le:a-cofib}, a map is a local-fibration
if and only if it has the local right lifting property with respect
to all maps $K \otimes X \inc L \otimes X$ for every
acyclic cofibration $K \inc L$ between finite simplicial sets.

\subsection{Sheaves of homotopy groups}

Let $F$ be a simplicial presheaf on
$\cC$.  Given an object $X$ of $\cC$ and a $0$-simplex $x$ in $F(X)$, we
define presheaves $\pi_n(F,x)$ on the site $\cC\ovcat X$
by the formula  
$U\mapsto \pi_n(F(U),x|_U)$.

\begin{defn}
\label{de:Ill-we}
A map of simplicial presheaves $f\colon F\ra G$ is a \dfn{local
weak equivalence} if
\begin{enumerate}[(1)]
\item The induced map $\pi_0 F \ra \pi_0 G$
yields an isomorphism upon sheafification, and
\item For every $X$ in $\cC$ and every basepoint $x$ in $F_0(X)$, the
map of presheaves on $\cC\ovcat X$ given by $\pi_n(\Ex^\infty F,x)\ra
\pi_n(\Ex^\infty G,fx)$ also becomes an isomorphism upon
sheafification.  (Here $\Ex^\infty F$ is the presheaf $U\mapsto
\Ex^\infty(F(U))$, of course).
\end{enumerate}
\end{defn}

Local weak equivalences are called `topological weak equivalences' in
\cite{J1}.  One can also use the presheaf $\pi_n^{loc}(F, x)$, whose
value on an object $U \map X$ is the set of based maps $S^n \map F(U)$
modulo the equivalence relation generated by local simplicial homotopy
(see \cite[p. 44]{J1}).  Two maps $S^n \map F(U)$ are locally
simplicially homotopic if there exists a covering sieve of $U$ such
that for every $V \map U$ in the sieve, the two restrictions $S^n \map
F(V)$ are simplicially homotopic as based maps.  The following result
appears in \cite[Prop.~1.18]{J1}, except for an unnecessary
hypothesis.

\begin{lemma}
\label{le:loc-sheaf}
The map $\pi_n (F,x) \map \pi_n^{loc} (F,x)$ is an isomorphism after
sheafification, for any simplicial presheaf $F$.
\end{lemma}

Before proving Lemma \ref{le:loc-sheaf}, we recall the following
property of sheafifications.

\begin{lemma}
\label{le:sheaf-iso}
A map $f\colon F\ra G$ between presheaves of sets
induces an isomorphism on sheafifications if and only if the following
two conditions are satisfied:
\begin{enumerate}[(1)]
\item Given any $X$ in $\cC$ and any $s$ in $G(X)$, there is a covering
sieve $R$ of $X$ such that 
the restriction $s\restr{U}$ belongs to the image of $F(U)$ in $G(U)$
for any element $U\ra X$ of $R$;
\item Given any $X$ in $\cC$ and any two sections $s$ and $t$ in $F(X)$ such
that $f(s)=f(t)$, there exists a covering sieve $R$ of $X$ such that
 $s\restr{U}=t\restr{U}$ in $F(U)$ for every element $U\ra X$ of $R$.
\end{enumerate}
\end{lemma}

\begin{proof}
Condition (2) is equivalent to $F^+\ra G^+$ being an objectwise
monomorphism, which in turn is equivalent to the same property for
$F^{++}\ra G^{++}$.  If $F^{++}\ra G^{++}$ is an objectwise surjection
then property (1) is easily seen to hold.  Finally, properties (1) and
(2) together imply that $\im(G(X)\ra G^{++}(X)) \subseteq
\im(F^{++}(X) \inc G^{++}(X))$.  From this one deduces that
$F^{++}\ra G^{++}$ is an objectwise surjection (using that the domain
and codomain are sheaves).
\end{proof}

We will make use of the above two conditions in studying sheaves of
homotopy groups.

\begin{proof}[Proof of Lemma \ref{le:loc-sheaf}]
Since local simplicial homotopy is a larger equivalence relation than
simplicial homotopy, the map is an objectwise surjection.  This verifies
condition (1) of Lemma \ref{le:sheaf-iso}.

For condition (2), suppose that $s$ and $t$ are two maps
$S^n \map F(U)$ that are related by a finite chain 
of local simplicially homotopies.
There is a finite sequence $s = s_0, s_1, \ldots, s_n = t$ of maps
$S^n \map F(U)$ such that $s_i$ and $s_{i+1}$ are simplicially homotopic
after restricting to a sieve $R_i$.  Taking $R$ to be a common refinement
of each $R_i$, we conclude that $s$ and $t$ are 
related by a chain of simplicial homotopies after restricting to $R$.
This verifies condition (2) of Lemma \ref{le:sheaf-iso}.
\end{proof}

In general, $\pi_n^{loc}(F,x)$ does not carry homotopical information
unless $F$ is locally-fibrant.  At first glance, the sheafification of
$\pi_n(F,x)$ seems not to be homotopically meaningful unless $F$ is
{\it objectwise\/} fibrant, but Lemma \ref{le:loc-sheaf} shows that we
only need $F$ to be locally-fibrant.  In other words, for
locally-fibrant simplicial presheaves one can ignore the presence of
$\Ex^\infty$ in Definition \ref{de:Ill-we}:

\begin{prop}
\label{pr:loc-fib-we}
If $F$ and $G$ are locally-fibrant, then a map $f:F\ra G$ is a local weak
equivalence if and only if 
\begin{enumerate}[(1)]
\item The induced map $\pi_0 F \ra \pi_0 G$
yields an isomorphism upon sheafification, and
\item For every $X$ in $\cC$ and every basepoint $x$ in $F_0(X)$, 
the map $\pi_n(F,x)\ra \pi_n(G,fx)$ is an
isomorphism upon sheafification.  
\end{enumerate}
\end{prop}

\begin{proof}
Consider the square
\[
\xymatrix{
F \ar[r] \ar[d] & \Ex^\infty F \ar[d] \\
G \ar[r] & \Ex^\infty G. }
\]
By \cite[Prop.~1.17]{J1} and Lemma~\ref{le:sheaf-iso}, the horizontal
maps satisfy the above conditions.  Therefore the left vertical map
satisfies the conditions if and only if the right vertical map does.
The usual complications with choosing basepoints do not arise because
$\Ex^\infty$ preserves $0$-simplices.
\end{proof}

\subsection{Local relative-homotopy-liftings}

The relative-homotopy-lifting criterion for weak equivalences of
simplicial sets (given in Proposition \ref{pr:ss-we}) has an obvious
extension to the presheaf category in which we only require local
liftings.

\begin{defn}
\label{de:llp}
Let $K\ra L$ be a map of simplicial sets.  A map $f\colon F\ra G$ of
simplicial presheaves is said to have the \mdfn{local RHLP} with
respect to $K\ra L$ if every diagram
\begin{myequation}
 \xymatrixcolsep{1.5pc}\xymatrix{
  &K \tens X \ar[r]\ar[d]\ar[dl] & F \ar[d] \\
  L\tens X \ar[d] \ar@{.>}[urr]  
             & L \tens X \ar[r]\ar[dl] & G  \\
  RH(L,K) \tens X  \ar@{.>}[urr]
}
\label{eq:loclift}
\end{myequation}
admits local liftings.
\end{defn}

In other words, the definition requires that 
there exists a covering sieve of $X$ such that for any
map $U \map X$ in the sieve, the induced diagram
\[
\xymatrix{
K \tens U \ar[r] \ar@{ >->}[d] & K \tens X \ar[r] & F \ar[d] \\
L \tens U \ar[r] \ar@{.>}[urr] & L \tens X \ar[r] & G       }
\]
has a relative-homotopy-lifting.  The liftings and simplicial
homotopies one gets as $U$ varies need not be compatible in any way.

The basic results about relative-homotopy-liftings from
Section~\ref{se:genprops} all go through in the present context.  One
only has to observe that the arguments require finitely many uses
of the lifting conditions.  The following result will be especially
useful to us:

\begin{lemma}
\label{le:lhlprops}
Let $f\colon F \ra G$ be a fixed map of simplicial presheaves.
Consider the class of all maps $K \map L$ of simplicial sets with
respect to which $f$ has the local RHLP.  This class is closed under
cobase changes and retracts.  If $G$ is locally-fibrant, then the
cofibrations in this class are also closed under composition.
\end{lemma}

\begin{proof}
The proof is the same as that of Lemma \ref{le:hlift-comp}, except
that the relative-homotopy-liftings are replaced by local
relative-homotopy-liftings.
\end{proof}

\begin{cor}
\label{co:lhlprops}
If $f\colon F\ra G$ has the local RHLP with respect to the maps 
$\bd{n}\inc \del{n}$ for all $n\geq 0$ and $G$ is locally-fibrant, 
then $f$ also has the local RHLP with respect to
all cofibrations $K\inc L$ of finite simplicial sets.
\end{cor}

\begin{proof}
Every cofibration $K \inc L$ can be constructed by a finite number
of compositions and cobase changes from the generating
cofibrations $\bd{n} \inc \del{n}$.
\end{proof}

Here is the main theorem of the paper:

\begin{thm}\mbox{}\par
\label{th:main}
\begin{enumerate}[(a)]
\item If $F$ and $G$ are locally-fibrant, then a map $F\ra G$ is a
local weak equivalence if and only if it has the local RHLP with
respect to the maps $\bd{n} \inc \del{n}$.  
\item If $F$ and $G$ are arbitrary and $\cR$ is any
fibrant-replacement functor for $\sSet$, then a map $F\ra G$ is a
local weak equivalence if and only if $\cR F \ra \cR G$ has the local
RHLP with respect to the maps $\bd{n} \inc \del{n}$.  
\end{enumerate}
\end{thm}

\begin{proof}
For (a), we begin by assuming that $F$ and $G$ are locally-fibrant and
that the map $f\colon F \ra G$ has the local RHLP.  By Corollary
\ref{co:lhlprops}, it has the local RHLP with respect to all
cofibrations between finite simplicial sets.  In particular, it has
the local RHLP with respect to $* \inc S^n$; this proves condition (1)
of Lemma \ref{le:sheaf-iso} for $\pi_n^{loc} (F, x) \map \pi_n^{loc}
(G, fx)$ (for $n=0$ one uses the RHLP with respect to $\emptyset\ra
*$).  On the other hand, $f$ also has the local RHLP with respect to
$S^n \Wedge S^n \inc RH(S^n, *)$; this proves condition (2) of Lemma
\ref{le:sheaf-iso} for $\pi_n^{loc} (F, x) \map \pi_n^{loc} (G, fx)$.
Here we use that local simplicial homotopy is an equivalence relation
for locally-fibrant simplicial presheaves \cite[Lem.~1.9]{J1}.  Now
Lemma \ref{le:loc-sheaf} and Proposition \ref{pr:loc-fib-we} tell us
that we have a local weak equivalence.

We now assume that $f\colon F\ra G$ is a local weak equivalence.
To prove that it has the local RHLP with respect to the maps $\bd{n}
\inc \del{n}$ we follow exactly the same argument as in Section
\ref{se:combprf}, observing that there are only finitely many
applications of the various lifting properties.

To prove (b), let $\cR$ be any fibrant-replacement functor for
$\sSet$.  We need only observe that $F\ra G$ is a local weak
equivalence if and only if $\cR F \ra \cR G$ is one.  Since $\cR F$
and $\cR G$ are objectwise-fibrant (hence locally-fibrant as well),
part (a) applies.
\end{proof}

%%%%%%%%%%%%%%%%%%%%%%%%%%%%%%%%%%%%%%%%%%%%%%%%%%%%%%%%%%%%%%%%%%%%%%

\section{Applications}

\label{se:app}

Both the lifting characterization of local weak equivalences and the
definition involving sheaves of homotopy groups are useful to have
around.  For instance, 
up to some technical difficulties in choosing basepoints,
it is transparent from the homotopy group
definition that the local weak equivalences have the two-out-of-three
property.  This is awkward to show using the lifting characterization,
however.  We now give some results which are easy consequences of our
lifting criterion.

\begin{prop}[Local right properness]
Let $F\ra G$ be a local weak equivalence between
simplicial presheaves, and let $J\ra G$ be a local-fibration.  Then
the map $J\times_G F \ra J$ is also a local weak equivalence.
\end{prop}

\begin{proof}
The proof is the same as the proof of Corollary \ref{co:ss-rp},
except that liftings are replaced by local liftings.
Note that $\Ex^\infty$ commutes with fibre-products of simplicial presheaves
since both $\Ex^\infty$ and fibre-products are defined objectwise.
Also, similar to observations in the proof of \cite[Prop.~1.17]{J1}, 
$\Ex^\infty$ preserves local-fibrations.
\end{proof}

Recall that a \dfn{local acyclic fibration} of simplicial presheaves
is a map that is both a local weak equivalence and a local-fibration.
The $(\Longrightarrow)$ direction of the following proposition was
proved in \cite[Lemma 7]{J3}---we can now prove the other one (see
\cite[Lemma 11]{J3} for a weaker version).

\begin{prop}
\label{pr:trlfib}
A map $p\colon F \map G$ of simplicial presheaves admits local liftings in 
every square
\begin{myequation}
\label{di:trfibsq}
\xymatrix{ \bd{n} \tens X \ar[r]\ar@{ >->}[d] & F \ar[d] \\
              \del{n} \tens X \ar@{.>}[ur]\ar[r] & G
}
\end{myequation}
if and only if it is a local acyclic fibration.
\end{prop}

\begin{proof}
First suppose that the local liftings exist.  Then $p$ also has the
local-lifting property with respect to the maps $\Lambda^{n,k} \inc
\Delta^n$, since these can be built from the maps $\bd{r}\ra\del{r}$
by finitely many cobase-changes and compositions.  Therefore
$p$ is a local-fibration.  Similar to observations in the proof of
\cite[Prop.~1.17]{J1}, the map $\Ex^\infty p\colon \Ex^\infty F \map
\Ex^\infty G$ has the local lifting property with respect to all maps
$\bd{n} \tens X \inc \Delta^{n} \tens X$.  In particular, $\Ex^\infty
p$ has the local RHLP with respect to the maps $\bd{n} \inc \del{n}$;
we use constant relative-homotopies.  Using Theorem~\ref{th:main}(b),
$p$ is a local weak equivalence.  This finishes one impliciation.

For the other direction, first assume that $F$ and $G$ are locally-fibrant.
Since $F\ra G$ is a local weak equivalence one gets local {\it
relative\/}-homotopy-liftings by Theorem~\ref{th:main}.  
Similar to the proof of Corollary \ref{co:ss-rp},
the fact that $F\ra G$ is a local fibration allows one 
to homotope the local homotopy-liftings to 
get actual local liftings.

Now suppose that $p$ is an arbitrary local acyclic fibration,
and suppose given a
lifting square as in (\ref{di:trfibsq}).  As we have already observed,
$\Ex^\infty p$ is also a local acyclic fibration, but with
locally-fibrant domain and codomain.  So by the previous paragraph there are
local liftings for the composite square
\[ \xymatrix{
\bd{n}\tens X \ar[r] \ar[d] & F \ar[r]\ar[d] & \Ex^\infty F \ar[d] \\
\del{n}\tens X \ar[r]\ar@{.>}[urr] & G \ar[r] & \Ex^\infty G.}
\]
This translates to saying that for a sufficiently large $k$ there are
local liftings in the square
\[\xymatrix{
\sd^k\bd{n} \tens X \ar[d]\ar[r]  
         &\bd{n}\tens X \ar[r]\ar[d] & F\ar[d] \\
\sd^k\del{n} \tens X \ar[r]\ar@{.>}[urr]  
         &\del{n}\tens X \ar[r] & G,}
\]
where the left horizontal maps are the `last vertex maps' \cite[p.~183]{GJ}.

Let $C$ be the mapping cylinder of $\sd^k\bd{n} \ra \bd{n}$, and let
$D$ be the mapping cylinder of $\sd^k\del{n} \ra \del{n}$.  Notice
that $C$ is a subcomplex of $D$.  Since the map $\sd^k\bd{n}\tens X
\ra F$ factors through $\bd{n}\tens X$, the constant homotopy $(
\sd^k\bd{n}\times\del{1})\tens X \ra F$ factors through $C\tens X$.
Now, we have squares
\[ \xymatrix{
[C \cup (\sd^k\del{n} \times\{0\})] \tens U \ar[r]\ar@{ >->}[d] 
      & F\ar[d] \\
D\tens U \ar[r] & G}
\]
for all $U\ra X$ in a covering sieve of $X$, where the maps $C\tens U
\ra F$ and $D\tens U\ra G$ are these `constant homotopies'.  The 
map $C \cup (\sd^k\del{n} \times\{0\}) \map D$
is a trivial cofibration
between finite simplicial sets (both the domain and codomain are
contractible), so the square has a local lifting.  By precomposing
these liftings with the inclusion $\del{n}\inc D$, one obtains local
liftings for the original square.
\end{proof}

\begin{cor}[cf. {\cite[Lemma 19]{J3}}]
\label{co:fsm7}
Let $p\colon F\ra G$ be a local-fibration, and let $i\colon K\inc L$
be a cofibration of finite simplicial sets.  If $p$ is a local weak
equivalence or $i$ is a weak equivalence, then the induced map
\[    F^L \ra F^K \times_{G^K} G^L 
\]
is a local acyclic fibration.
\end{cor}

\begin{proof}
To see that the map is a local acyclic fibration, it is enough by
Proposition~\ref{pr:trlfib} to check that it has the local lifting
property with respect to the maps $\bd{n}\ra \del{n}$.  By
adjointness, one need only check that $F\ra G$ has the local lifting
property with respect to the map
\[ j\colon(L\times\bd{n}) \cup (K\times \del{n}) \ra L\times \del{n}.\]
If $K\ra L$ is an acyclic cofibration then so is $j$, and therefore
the result follows from Lemma~\ref{le:a-cofib} and the definition of
local fibration.  If $F\ra G$ was an acyclic fibration, then the
result follows from Proposition~\ref{pr:trlfib} because $j$ is
obtained by a finite number of cobase changes and compositions from
the inclusions $\bd{n}\ra\del{n}$.
\end{proof}

Finally, we end with the following result.  It is needed in
\cite{DHI}, and we know of no proof that avoids the local lifting
techniques we've just developed.

\begin{cor}
Let $F\ra G$ be a map between locally-fibrant simplicial presheaves,
and let $K\inc L$ be a cofibration of finite simplicial sets.  If
either map is a weak equivalence then the induced map $F^L \ra F^K
\times_{G^K} G^L$ is a weak equivalence.
\end{cor}

\begin{proof}
First, we know from \cite[Cor. 1.5]{J1} that both $F^L$ and
$F^K\times_{G^K} G^L$ are locally-fibrant.  A lifting square
\[\xymatrix{\bd{n} \tens X \ar[r]\ar@{ >->}[d] & F^L \ar[d] \\
             \del{n} \tens X \ar[r] & F^K\times_{G^K} G^L}
\] 
may be rewritten via adjointness as
\[\xymatrix{ 
M \tens X \ar[r]\ar@{ >->}[d] & F \ar[d] \\
(L\times \del{n}) \tens X \ar[r] & G,}
\]
where $M=(L\times \bd{n}) \cup (K\times\del{n})$.  The map $M\inc
L\times\del{n}$ is a cofibration between finite simplicial sets.  When
$F \map G$ is a local weak equivalence, Corollary \ref{co:lhlprops}
and Theorem~\ref{th:main}(a) tell us that the above square has a
relative-homotopy-lifting.  Using adjointness once again, we get a
relative-homotopy-lifting for the original square.

The other case is similar.  If $K\inc L$ is a weak equivalence then
$M\inc L\times\del{n}$ is also one.  So by the local version of
Proposition~\ref{pr:fibhl} (which only works for acyclic cofibrations
between finite simplicial sets), we have a relative-homotopy-lifting
since $F$ and $G$ are locally-fibrant.
\end{proof}

%%%%%%%%%%%%%%%%%%%%%%%%%%%%%%%%%%%%%
\bibliographystyle{amsalpha}

\end{document}